\documentclass{conm-p-l}

\newtheorem{theorem}{Theorem}[section]
\newtheorem{lemma}[theorem]{Lemma}
\newtheorem{prop}[theorem]{Proposition}
\newtheorem{coro}[theorem]{Corollary}

\theoremstyle{definition}
\newtheorem{definition}[theorem]{Definition}
\newtheorem{nota}[theorem]{Notation}
\newtheorem{example}[theorem]{Example}

\theoremstyle{remark}

\def\ZZ{\mathbb{Z}}

\def\calG{\mathcal{G}}
\def\calH{\mathcal{H}}

\def\NN{\mathbb{N}}
\def\ZZ{\mathbb{Z}}

\def\RR{\mathbb{R}}

\def\span{\overline{\text{\rm span}}}
\def\ospan{{\text{\rm span}}}
\def\NN{\mathbb{N}}
\def\calH{\mathcal{H}}
\def\H{\mathcal{H}}

\newcommand{\ip}[2]{\left\langle#1,#2\right\rangle}
\newcommand{\absip}[2]{\left| \left\langle#1,#2\right\rangle \right|}
\newcommand{\norm}[1]{\left\lVert#1\right\rVert}

\numberwithin{equation}{section}

\begin{document}

\title{Frames of subspaces}

\author{Peter G. Casazza}
\address{Department of Mathematics,
University of Missouri,
Columbia, Missouri 65211 USA}
\email{pete@math.missouri.edu}
\thanks{The first author was supported by NSF DMS 0102686.}

\author{Gitta Kutyniok}
\address{Institute of Mathematics,
University of Paderborn,
33095 Paderborn, Germany}
\email{gittak@uni-paderborn.de}
\thanks{The second author was supported by Forschungspreis 2003 der
Universit\"at Paderborn.}

\subjclass{Primary 42C15; Secondary 46C99}
\date{October 21, 2004.}


\keywords{Abstract frame theory, frame, harmonic frame,
Hilbert space, resolution of the identity, Riesz basis, Riesz frame}

\begin{abstract}
One approach to ease the construction of frames is to first construct local
components and then build a global frame from these.
In this paper we will show that the study of the relation between a frame and its
local components leads to the definition of a frame of subspaces. We introduce this
new notion and prove that it provides us with the link we need.
It will also turn out that frames of subspaces behave as a generalization of frames.
In particular, we can
define an analysis, a synthesis and a frame operator for a frame of subspaces, which
even yield a reconstruction formula. Also concepts such as completeness, minimality, and
exactness are introduced and investigated. We further study several constructions of
frames of subspaces, and also of frames and Riesz frames using the theory of frames
of subspaces.
An important special case are harmonic frames of subspaces which generalize harmonic
frames. We show that wavelet subspaces coming from multiresolution analysis belong to this class.
\end{abstract}

\maketitle


\section{Introduction}

During the last 20 years the theory of frames has been growing rapidly, since several
new applications have been developed.
For example, besides traditional applications as signal processing, image
processing, data compression, and sampling theory, frames are now used
to mitigate the effect of losses in packet-based communication
systems and hence to improve the robustness of data transmission
\cite{CK01,GKK01}, and to design high-rate constellations with full
diversity in multiple-antenna code design \cite{HHSS01}.

To handle these emerging applications of frames new methods have to be developed.
One starting point is to first build frames ``locally'' and then piece them together
to obtain frames for the whole space. One advantage of this idea is that it would
facilitate the construction of frames for special applications, since we can first
construct frames or choose already known frames for smaller spaces. And in a second
step one would construct a frame for the whole space from them.
Therefore it is necessary to derive conditions for these components, so that there exists
a construction, which yields a frame for the whole space with special properties.
Various approaches to piecing together familes of vectors to get a frame for the
whole space have been done over the years going back to Duffin and Schaeffer's
original work \cite{DS}.  One approach used in the wavelet as well as in the Gabor case \cite{DGM,ACM03}
is to start with non-frame sequences and piece them together to build frames for the whole
space.  Another is to build frames locally and piece them together orthogonally to get
frames. We refer to Heil and Walnut \cite{HW89} for an excellent introduction to these methods
and Gabor frames in general.  Recently, another approach was introduced by
Fornasier \cite{F,F1}.
Fornasier uses subspaces which are quasi-orthogonal to construct local frames and piece them
together to get global frames.

In this paper we will formulate a general method for piecing together local frames to get
global frames.  The importance of this approach is that it is both necessary and
sufficient for the the construction of global frames from local frames. Some of these
results are generalizations of Fornasier's work \cite{F,F1} although they were done
before his papers became available to us.

Another motivation comes from the theory of C*-algebras. Just recently Ca\-sazza, Christensen, Lindner,
and Vershynin \cite{CCLV03} proved that the so-called "Feichtinger conjecture" is equivalent
to the weak Bourgain-Tzafriri conjecture. The Feichtinger conjecture states that each bounded
frame is a finite union of Riesz basic sequences.  Then, Ca\-sazza and Vershynin \cite{CV03}
showed that the Kadison-Singer problem is equivalent to the strong Bourgain-Tzafriri
conjecture and that these two problems have a positive solution if and only if both
the Feichtinger conjecture and the $F_{\epsilon}$-conjecture have positive solutions.
The $F_{\epsilon}$-conjecture states:  For every $\epsilon >0$, every unit norm Riesz basis is
a finite union of $(1+\epsilon)$-unconditional basic sequences.  A unit norm sequence $\{f_i\}_{i\in I}$
is a $(1+\epsilon)$-basic sequence if for every sequence of scalars $\{a_i \}_{i\in I}$
we have
$$
(1-\epsilon )\sum_{i\in I}|a_i |^2 \le \left\|\sum_{i\in I}a_i f_i \right\|^2 \le
(1+\epsilon)
\sum_{i\in I}|a_i |^2 .
$$
To attack these problems it is important to know into which components we can divide
a frame.  As we will see in this paper, the necessary divisions will form a frame of
subspaces for the space.  At this time, it is not even known how to divide a frame into
two infinite frame sequences.

In this paper we want to answer the following two questions, which relate to the two different
motivations:
\begin{itemize}
\item Let $\{W_i\}_{i \in I}$ be a collection of closed subspaces in a Hilbert space $\H$
in which we want to decompose our function, where each subspace $W_i$ is equipped
with a weight $v_i$, which indicate its importance.
When can we find frames for $W_i$ for each $i \in I$ so that the
collection of all of them is a frame with special properties for the whole space $\H$?
\item Let $\{f_i\}_{i \in I}$ be a frame for a Hilbert space $\H$, and let
$I=\bigcup_{j \in \ZZ} I_j$ be a partition of $I$ so that $\{f_i\}_{i \in I_j}$
is a frame sequence for each $j \in \ZZ$. Which relations exist between
the closed linear spans of $\{f_i\}_{i \in I_j}$, $j \in \ZZ$ ?
\end{itemize}

\medskip

We start our consideration by giving a brief review of the definitions and basic
properties of frames and bases and stating some notation in Section \ref{review}.

In Section \ref{fos} it will turn out that both questions above lead to the definition
of a frame of subspaces. In the first subsection we will state the definition of
a frame of subspaces for a given family of closed subspaces $\{W_i\}_{i \in I}$
in a Hilbert space and a family of weights $\{v_i\}_{i \in I}$. Then it is shown that this
definition leads to some answers to the above questions (see Theorem \ref{transfer_frame}),
since it shows that frames of subspaces behave as a link between local components of a frame
and the global structure.
This will also enlighten the advantage of our approach, since now we can choose the frames
for the single subspaces $W_i$ arbitrarily and always get a frame for the whole Hilbert
space by just collecting them together. Thus it differs from previous
approaches and is a generalization of the approach of Fornasier \cite{F,F1}.
It will turn out that frames of subspaces behave as a generalization
of frames. We first give a definition of completeness of a family of subspaces and show that the
relation between this property and the notion of a frame of subspaces
is similar to the relation between the definition of completeness of a sequence and
a frame.
Further in Subsection \ref{frame_prop} we introduce an analysis and a synthesis operator,
a frame operator, and a dual frame of subspaces for a given frame of subspaces and prove
that they behave in an analogous way as the corresponding objects in abstract frame theory.
We even obtain a reconstruction formula using these ingredients (Proposition \ref{reconstruction}).
The next subsection deals with Parseval frames of subspaces, which share several
properties with Parseval frames. Finally in Subsection \ref{roti} we show that using
the theory of frames of subspaces we can construct several useful resolutions of the identity.

Section \ref{Riesz_dec} deals with Riesz decompositions, which are a generalization of
the notion of Riesz bases to our general setting. We further define minimality for
a family of subspaces and show that it behaves as expected. Also exactness
is defined in a canonical way. However, it will turn out that this property is
much weaker than exactness of a frame (compare Theorem \ref{Riesz_minimal_exact}).

Some constructions are given in Section \ref{constr}. Here we first state some
results which help constructing frames of subspaces. An extended example concerning
the situation of Gabor frames is added.
In Subsection \ref{constr_frame_Riesz} we then show how to construct frames and Riesz
frames using a frame of subspaces.

Finally, Section \ref{harmonic_fos} deals with harmonic frames of subspaces. These
are a generalization of harmonic frames, which distinguish themselves
by having an easy construction formula. In both
the finite and the infinite dimensional cases we give the definition of a harmonic frame of subspaces,
state some results, and give examples, e.g., subspaces coming from Gabor systems
and subspaces coming from multiresolution analysis, for their occurance.


\section{Review of frames and some notation}
\label{review}

First we will briefly recall the definitions and basic properties
of frames and bases. For more information we refer to the survey
articles by Casazza \cite{Cas01,Cas00}, the books by Christensen
\cite{Chr03}, Gr\"ochenig \cite{Gro01}, and Young \cite{You80} and the
research-tutorial by Heil and Walnut \cite{HW89}.

Let $\calH$ be a separable Hilbert space and let $I$ be an indexing set.
A family $\{f_i\}_{i \in I}$ is a {\em frame} for $\calH$, if there
exist $0 < A \le B < \infty$ such that for all $h \in \calH$,
\begin{equation}
\label{framedef}
A \, \norm{h}^2 \le \sum_{i \in I} \absip{h}{f_i}^2
\le B \, \norm{h}^2.
\end{equation}
The constants $A$ and $B$ are called a {\em lower} and {\em upper frame bound}
for the frame.
Those sequences which satisfy only the upper inequality in (\ref{framedef})
are called {\em Bessel sequences}.
A frame is {\em tight}, if $A=B$. If $A=B=1$, it is called a
{\em Parseval frame}. We call a frame $\{f_i\}_{i \in I}$ {\em uniform}
(or {\em equal norm}), if
we have $\|f_i\| = \|f_j\|$ for all $i,j \in I$.
A frame is {\em exact}, if it ceases to be a frame
whenever any single element is deleted from the sequence $\{f_i\}_{i \in I}$.
We say that two frames $\{f_i\}_{i \in I}$, $\{g_i\}_{i \in I}$ for $\calH$ are
{\em equivalent}, if there exists an invertible operator $U : \calH \to \calH$
satisfying $Uf_i = g_i$ for all $i \in I$. If $U$ is a unitary operator,
$\{f_i\}_{i \in I}$ and $\{g_i\}_{i \in I}$ are called {\em unitarily equivalent}.
The {\em synthesis operator} $T_f : l^2(I) \to \H$ of a frame $f=\{f_i\}_{i \in I}$ is
defined by $T_f(c) = \sum_{i \in I} c_i f_i$ for each sequence of scalars
$c=\{c_i\}_{i \in I} \in l^2(I)$. The adjoint operator $T^*_f : \H \to l^2(I)$,
the so-called {\em  analysis operator} of $f=\{f_i\}_{i \in I}$,
is given by $T^*_f(g)=\{\ip{f_i}{g}\}_{i \in I}$.
Then the {\em frame operator} $S_f (h) = T_f T_f^*(h) = \sum_{i \in I} \ip{h}{f_i}f_i$
associated with $\{f_i\}_{i \in I}$ is a bounded, invertible, and positive operator
mapping $\calH$ onto itself. This provides the reconstruction formula
\[h = S_f^{-1}S_f(h) = \sum_{i \in I} \ip{h}{f_i} \tilde{f_i} = \sum_{i \in I}
\langle h,\tilde{f_i}\rangle f_i,\]
where $\tilde{f_i} = S_f^{-1}f_i$. The family $\{\tilde{f_i}\}_{i \in I}$ is also
a frame for $\calH$, called the {\em canonical dual frame} of $\{f_i\}_{i \in I}$.
A sequence is called a {\em frame sequence}, if it is a frame only for its closed
linear span. Moreover, we say that a frame $\{f_i\}_{i \in I}$ is a {\em Riesz frame},
if every subfamily of the sequence $\{f_i\}_{i \in I}$ is a frame sequence with uniform
frame bounds $A$ and $B$.

As important example of frames are the so-called {\em harmonic frames,} which are
uniform Parseval frames of the form $\{U^i \varphi\}_{i \in I}$, where $U$ is a
unitary operator on $\H$ and $I = \{0,\ldots,N-1\}$, $N \in \NN$ or $I = \ZZ$.
Concerning a classification of harmonic frames we refer to the paper by Casazza and Kova\^cevi\'c
\cite{CK01}.

{\em Riesz bases} are special cases of frames, and can be characterized as those frames
which are biorthogonal to their dual frames. An equivalent definition is the following.
A family $\{f_i\}_{i \in I}$ is a Riesz basis for $\calH$, if there exist $0 < A \le B
< \infty$ such that for all sequences of scalars $c = \{c_i\}_{i \in I}$,
\[ A \, \norm{c}_2 \le \left\|\sum_{i \in I} c_i f_i\right\| \le B \, \norm{c}_2.\]
We define the {\em Riesz basis constants} for $\{f_i\}_{i \in I}$ to be the largest
number $A$ and the smallest number $B$ such that this inequality holds for all
sequences of scalars $c$. If $\{f_i\}_{i \in I}$ is a Riesz basis
only for its closed linear span, we call it a {\em Riesz basic sequence}.

An arbitrary sequence $\{f_i\}_{i \in I}$ in $\calH$ is {\em minimal}, if
$f_i \not\in \span_{j \in I, j\not= i}\{f_{j}\}$ for all $i \in I$, or equivalently
if there exists a sequence $\{\tilde{f_i}\}_{i \in I}$, which is biorthogonal to
$\{f_i\}_{i \in I}$. It is
{\em complete}, if the span of $\{f_i\}_{i \in I}$ is dense in $\calH$.

\medskip

We conclude this section by giving some notation and remarks. Throughout
this paper $\H$ shall always denote an arbitrary separable Hilbert space.
Furthermore all subspaces are assumed to be closed although this is not stated
explicitely. Moreover, for the remainder a sequence $\{v_i\}_{i \in I}$ always
denotes a family of weights, i.e., $v_i > 0$ for all $i \in I$.

In addition we use the following notation. Dependent on the context $I$
denotes an indexing set or the identity operator.
If $W$ is a subspace of a Hilbert space $\H$,
we let ${\pi}_{W}$ denote the orthogonal projection
of $\H$ onto $W$.  If $\{e_i \}_{i\in I}$ is an orthonormal
basis for $\H$ and $J\subset I$, ${\pi}_{J}$ is the orthogonal
projection of $\H$ onto $\span_{i \in J} \{e_{i}\}$.


\section{Frames of subspaces}
\label{fos}

\subsection{Definition and basic properties}
\label{def_basic_prop}

We start with the definition of a frame of subspaces. It will turn out that frames
of subspaces share many of the properties of frames, and thus can be viewed as a
generalization of frames.

\begin{definition}
Let $I$ be some index set, and let $\{v_i\}_{i \in I}$ be a family
of weights, i.e., $v_i > 0$ for all $i \in I$.
A family of closed subspaces $\{W_i\}_{i\in I}$ of a Hilbert space $\H$ is a
{\it frame of subspaces with respect to $\{v_i\}_{i \in I}$ for $\H$}, if there exist
constants $0 < C \le D < \infty$ such that
\begin{equation} \label{deffos}
C\|f\|^2 \le \sum_{i\in I} v_i^2 \|{\pi}_{W_i}(f)\|^2 \le D\|f\|^2
\ \ \mbox{for all $f\in \H$}.
\end{equation}
We call $C$ and $D$ the {\em frame bounds} for the frame of subspaces.
The family $\{W_i\}_{i\in I}$ is called a {\em $C$-tight frame of subspaces
with respect to $\{v_i\}_{i \in I}$}, if in (\ref{deffos})
the constants $C$ and $D$ can be chosen so that $C=D$, a {\em Parseval frame
of subspaces with respect to $\{v_i\}_{i \in I}$} provided that $C=D=1$
and an {\em orthonormal basis of subspaces} if $\H = \bigoplus_{i \in I} W_i$.
Moreover, we call a frame of subspaces with
respect to $\{v_i\}_{i \in I}$ {\em $v$-uniform}, if $v := v_i = v_j$ for all $i,j \in I$.
If we only have the upper bound, we call $\{W_i\}_{i\in I}$
a {\it Bessel sequence of subspaces with respect to $\{v_i\}_{i \in I}$} with {\it Bessel bound}
$D$.
\end{definition}

Condition (\ref{deffos}) states the necessary (and also sufficient) interaction
between the subspaces so that taking frames from them and putting them together
yields a frame for the whole space.

The importance of this definition is that it is both necessary and sufficient for us to
be able to string together frames for each of the subspaces $W_i$ (with uniformly bounded
frame constants) to get a frame for $\H$.  This is contained in the next theorem.
The implication (3) $\Rightarrow$ (1) of the following result is
\cite[Proposition 4.5]{CK01}.  Fornasier \cite{F,F1} obtains a similar result
for quasi-orthogonal decompositions.

\begin{theorem} \label{transfer_frame}
For each $i \in I$ let $v_i > 0$ and let $\{f_{ij}\}_{j \in J_i}$ be
a frame sequence in $\mathcal{H}$ with frame bounds $A_i$ and $B_i$. Define
$W_i = \span_{j\in J_i}\{f_{ij}\}$ for all $i\in I$ and
choose an orthonormal basis $\{e_{ij}\}_{j \in J_i}$ for each subspace $W_i$.
Suppose that $0 < A = \inf_{i \in I} A_i \le B = \sup_{i \in I} B_i < \infty$.
The following conditions are equivalent.
\begin{enumerate}
\item $\{v_i f_{ij}\}_{i \in I,j \in J_i}$ is a frame for $\mathcal{H}$.
\item $\{v_i e_{ij}\}_{i \in I,j \in J_i}$ is a frame for $\mathcal{H}$.
\item $\{W_{i}\}_{i\in I}$ is a frame of subspaces with respect to $\{v_i\}_{i \in I}$
for $\H$.
\end{enumerate}
\end{theorem}

\begin{proof}
Since for each $i \in I$, $\{f_{ij}\}_{j \in J_i}$ is a frame for $W_i$
with frame bounds $A_i$ and $B_i$, we obtain
\[ A \sum_{i \in I}v_i^2\norm{{\pi}_{W_i} (f)}^2
\le \sum_{i \in I} A_i v_i^2\norm{{\pi}_{W_i} (f)}^2
\le \sum_{i \in I} \sum_{j \in J_i} \absip{{\pi}_{W_i} (f)}{v_i f_{ij}}^2\]
\[\le \sum_{i \in I} B_i v_i^2 \norm{{\pi}_{W_i} (f)}^2
\le B \sum_{i \in I} v_i^2 \norm{{\pi}_{W_i} (f)}^2.\]
Now we observe that
\[\sum_{i \in I} \sum_{j \in J_i} \absip{{\pi}_{W_i} (f)}{v_i f_{ij}}^2
= \sum_{i \in I} \sum_{j \in J_i} \absip{f}{v_i f_{ij}}^2.\]
This shows that provided $\{v_i f_{ij}\}_{i \in I,j \in J_i}$ is a frame for $\mathcal{H}$
with frame bounds $C$ and $D$, the sets $\{W_i\}_{i \in I}$ form a frame
of subspaces with respect to $\{v_i\}_{i \in I}$ for $\calH$ with frame
bounds $\frac{C}{B}$ and $\frac{D}{A}$. Moreover, if $\{W_i\}_{i \in I}$ is a frame
of subspaces with respect to $\{v_i\}_{i \in I}$ for $\calH$ with frame
bounds $C$ and $D$, the calculation above implies that
$\{v_i f_{ij}\}_{i \in I,j \in J_i}$ is a frame for $\mathcal{H}$
with frame bounds $AC$ and $BD$. Thus (1) $\Leftrightarrow$ (3).

To prove the equivalence of (2) and (3), note that we can now
actually calculate the orthogonal projections in the following way
\[ v_i^2 \norm{{\pi}_{W_i} (f)}^2
= v_i^2 \left\| \sum_{j \in J_i} \ip{f}{e_{ij}}e_{ij} \right\|^2
= \sum_{j \in J_i} \absip{f}{v_i e_{ij}}^2.\]
From this the claim follows immediately.
\end{proof}

The definition of completeness of a sequence gives rise to a definition
of completeness for a sequence of subspaces.

\begin{definition}
A family of subspaces $\{W_i\}_{i\in I}$ of $\H$
is called {\em complete}, if
\[\span_{i \in I} \{W_i\} = \H.\]
\end{definition}

The next lemma possesses a well-known analog in the frame situation.

\begin{lemma}\label{fos-complete}
Let $\{W_i\}_{i\in I}$ be a family of subspaces in $\H$, and let
$\{v_i\}_{i \in I}$ be a family of weights.
If $\{W_i\}_{i\in I}$ is a frame of subspaces with respect to $\{v_i\}_{i \in I}$
for $\H$, then it is complete.
\end{lemma}

\begin{proof}
Assume that $\{W_i\}_{i\in I}$ is not complete. Then there exists some
$f \in \H$, $f \neq 0$ with $f \perp \span_{i \in I} \{W_i\}$. It follows that
$\sum_{i \in I}v_i^2 \|\pi_{W_i}(f)\|^2 = 0$, hence $\{W_i\}_{i\in I}$ is not
a frame of subspaces.
\end{proof}

To check completeness of a frame of subspaces, we derive the following useful
characterization.

\begin{lemma} \label{transfer_complete}
Let $\{W_i\}_{i\in I}$ be a family of subspaces in $\H$, and for each $i \in I$
let $\{e_{ij}\}_{j \in J_i}$ be an orthonormal basis for $W_i$.
Then the following conditions are equivalent.
\begin{enumerate}
\item $\{W_i\}_{i\in I}$ is complete.
\item $\{e_{ij}\}_{i \in I, j \in J_i}$ is complete.
\end{enumerate}
\end{lemma}

\begin{proof}
The equivalence of (1) and (2) follows immediately from the definitions.

\end{proof}

If we remove an element from a frame, we obtain either another frame or an incomplete
set \cite[Theorem 5.4.7]{Chr03}. A similar result holds in our situation.

\begin{prop}
The removal of a subspace from a frame of subspaces with respect to some family of weights
leaves either a frame of subspaces with respect to the same family of weights
or an incomplete family of subspaces.
\end{prop}

\begin{proof}
Let $\{W_i\}_{i\in I}$ be a frame of subspaces with respect to $\{v_i\}_{i \in I}$ for $\H$, and
for each $i \in I$ let $\{e_{ij}\}_{j \in J_i}$ be an orthonormal basis for $W_i$.
By Theorem \ref{transfer_frame}, $\{v_i e_{ij}\}_{i \in I,j \in J_i}$ is a frame for $\H$.
Let $i_0 \in I$. By \cite[Theorem 5.4.7]{Chr03}, $\{v_i e_{ij}\}_{i \in I \backslash \{i_0\},j \in J_i}$
is either a frame or an incomplete set. If it is a frame, again by Theorem \ref{transfer_frame},
also $\{W_i\}_{i\in I\backslash \{i_0\}}$ is a frame of subspaces with respect to $\{v_i\}_{i \in I}$
for $\H$. Now suppose that $\{v_i e_{ij}\}_{i \in I \backslash \{i_0\},j \in J_i}$ and
hence $\{e_{ij}\}_{i \in I \backslash \{i_0\},j \in J_i}$
is an incomplete set. By Lemma \ref{transfer_complete}, also $\{W_i\}_{i\in I\backslash \{i_0\}}$ is
incomplete.
\end{proof}

We further observe that the intersection of the elements of a frame of subspaces
with a subspace still leaves a frame of subspaces for a smaller space.

\begin{lemma}
Let $V$ be a subspace of $\calH$ and let $\{ W_i \}_{i\in I}$ be a frame of
subspaces with respect to $\{v_i\}_{i \in I}$ for $\calH$ with frame bounds
$C$ and $D$. Then $\{W_i \cap V\}_{i\in I}$ is a frame of
subspaces with respect to $\{v_i\}_{i \in I}$ for $V$ with frame bounds
$C$ and $D$.
\end{lemma}

\begin{proof}
For all $f \in V$ we have
\[ \sum_{i\in I} v_i^2 \|{\pi}_{W_i}(f)\|^2 = \sum_{i\in I} v_i^2 \|{\pi}_{W_i\cap V}(f)\|^2.\]
From this the result follows at once.
\end{proof}


\subsection{Frame properties}
\label{frame_prop}

In this subsection we will show that a frame of subspaces behaves as a
generalization of a frame, thus providing an associated analysis and synthesis
operator, a frame operator and a dual object.

For the definition of an analysis and a synthesis operator for a frame of
subspaces, we will need the following notation.

\begin{nota}
For each family of subspaces $\{ W_i \}_{i\in I}$ of $\H$, we define the
space $\left ( \sum_{i\in I} \oplus W_i \right ) _{{\ell}_{2}}$ by
\[\left ( \sum_{i\in I} \oplus W_i \right ) _{{\ell}_{2}} =
\{\{f_i \}_{i\in I}| f_i \in W_i  \mbox{ and } \sum_{i\in I}
\|f_i \|^2 < \infty \}\]
with inner product given by
\[\langle \{f_i \}_{i\in I} , \{g_i \}_{i\in I}\rangle =
\sum_{i\in I}\langle f_i, g_i \rangle.\]
\end{nota}

We start with the definition of a synthesis operator for a frame of subspaces. To show that
the series appearing in this formula converges unconditionally, we need the
next lemma.

\begin{lemma} \label{convergence}
Let $\{W_i\}_{i \in I}$ be a Bessel sequence of subspaces with respect to $\{v_i\}_{i \in I}$
for $\H$. Then, for each sequence $\{f_i\}_{i \in I}$ with $f_i \in W_i$ for each $i \in I$,
the series $\sum_{i \in I} v_i f_i$ converges unconditionally.
\end{lemma}

\begin{proof}
Let $f = \{f_i \}_{i\in I}\in \left ( \sum_{i \in I} \oplus W_i  \right)_{{\ell}_{2}}$.
Fix $J\subset I$ with $|J|< \infty$ and let $g = \sum_{i\in J}v_i f_i$.
Then we compute
\[\|\sum_{i\in J} v_i f_i \|^4 = \left( \langle g, \sum_{i\in J}v_i f_i \rangle
\right) ^2
= \left( \sum_{i\in J} v_i \langle  {\pi}_{W_i}(g), f_i \rangle \right) ^2
\le \left( \sum_{i\in J} v_i \|{\pi}_{W_i}(g) \|\|f_i \| \right) ^2\]
\[\le \sum_{i\in J} v_i^2 \|{\pi}_{W_i}(g) \|^2 \sum_{i\in J} \|f_i \|^2 \le
D\|g\|^2 \sum_{i\in J}\|f_i \|^2
\le D \|\sum_{i\in J} v_i f_i \|^2 \|f\|^2.\]
Hence,
\[\|\sum_{i\in J}v_i f_i \|^2 \le D\|f\|^2.\]
It follows that $\sum_{i\in I}v_i f_i $ is weakly unconditionally Cauchy and hence
unconditionally convergent in $\H$ (see \cite{D}, page 44, Theorems 6 and 8).
\end{proof}

\begin{definition}
Let $\{W_i \}_{i\in I}$ be a frame of subspaces with respect to $\{v_i\}_{i \in I}$ for
$\H$. Then the {\em synthesis operator} for $\{W_i \}_{i\in I}$ and $\{v_i\}_{i \in I}$
is the operator
\[T_{W,v}: \left ( \sum_{i\in I} \oplus W_i \right )_{{\ell}_{2}} \longrightarrow \H \]
defined by
\[T_{W,v}(f) = \sum_{i\in I}v_i f_i \quad \mbox{for all }
f = \{f_i \}_{i\in I} \in (\sum_{i\in I}\oplus W_i )_{{\ell}_{2}}.\]
We call the adjoint $T_{W,v}^{*}$ of the synthesis operator the {\it analysis operator}.
\end{definition}

The following proposition will provide us with a concrete formula for the analysis operator.

\begin{prop}
Let $\{W_i \}_{i\in I}$ be a frame of subspaces with respect to $\{v_i\}_{i \in I}$
for $\H$. Then the analysis operator $T_{W,v}^{*}:\H \rightarrow
(\sum_{i\in I}\oplus W_i )_{{\ell}_{2}}$ is given by
\[T_{W,v}^{*}(f) = \{v_i{\pi}_{W_i}(f)\}_{i\in I}.\]
\end{prop}

\begin{proof}
Let $f \in \H$ and $g = \{g_i \}_{i\in I} \in (\sum_{i\in I}\oplus W_i)_{{\ell}_{2}}$.
Using the definition of $T_{W,v}$ we compute that
\[\langle T_{W,v}^{*}(f) , g\rangle = \langle f, T_{W,v}(g) \rangle =
\langle f, \sum_{i\in I}v_i g_i \rangle
= \sum_{i\in I}v_i\langle f, g_i \rangle.\]
Since $g_i \in W_i $ for each $i \in I$, we can continue in the following way:
\[\sum_{i\in I}v_i\langle f, g_i \rangle
= \sum_{i\in I}v_i\langle \pi_{W_i}(f) , g_i \rangle
= \langle \{v_i{\pi}_{W_i}(f)\}_{i\in I}, \{g_i\}_{i \in I}\rangle.\]
\end{proof}

The well-known relations between a frame and the associated analysis and
synthesis operator also holds in our more general situation.

\begin{theorem}\label{PP1}
Let $\{W_i \}_{i\in I}$ be a family of subspaces in
$\H$, and let $\{v_i\}_{i \in I}$ be a family of weights.
Then the following conditions are equivalent.
\begin{enumerate}
\item  $\{W_i \}_{i\in I}$ is a frame of subspaces with respect to $\{v_i\}_{i \in I}$ for $\H$.
\item The synthesis operator $T_{W,v}$ is bounded, linear and onto.
\item  The analysis operator $T^*_{W,v}$ is a (possibly into) isomorphism.
\end{enumerate}
\end{theorem}

\begin{proof}
First we prove (1) $\Leftrightarrow$ (3). This claim follows immediately from the
fact that for each $f\in \H$ we have
\[\|T_{W,v}^*(f)\|^2 = \|\{v_i{\pi}_{W_i}(f)\}_{i\in I}\|^2 = \sum_{i\in I}
v_i^2\|{\pi}_{W_i}(f)\|^2.\]
Further recall that (2) $\Leftrightarrow$ (3) holds in general for
each operator on a Hilbert space.
\end{proof}

In an analogous way as in frame theory we can define equivalence classes of
frames of subspaces. Using the synthesis operator we can also characterize
exactly the elements belonging to the same equivalence class.

\begin{definition}
Let $\{W_i \}_{i\in I}$ and $\{\widetilde{W}_i \}_{i\in I}$ be frames of subspaces with respect to
the same family of weights. We say
that they are {\em (unitarily) equivalent}, if there exists an (unitary)
invertible operator $U$ on $\calH$ such that $W_i = U(\widetilde{W}_i)$ for all $i \in I$.
\end{definition}

\begin{lemma}
Let $\{W_i \}_{i\in I}$ and $\{\widetilde{W}_i \}_{i\in I}$ be frames of subspaces with respect to
the same family of weights $\{v_i\}_{i \in I}$.
The following conditions are equivalent.
\begin{enumerate}
\item $\{W_i \}_{i\in I}$ and $\{\widetilde{W}_i \}_{i\in I}$ are (unitarily) equivalent.
\item There exists an (unitary) invertible operator $U$ on $\calH$ such that
$T_{W,v} = U^{-1} T_{\widetilde{W},v} U$, where $U$ is applied to each component.
\end{enumerate}
\end{lemma}

\begin{proof}
This follows immediately from the definition of the synthesis operator.
\end{proof}

As in the well-known frame situation, there also exists an associated frame operator
for each frame of subspaces which satisfies similar properties as we will see in
the next proposition. For instance we even obtain a reconstruction formula.

\begin{definition}
Let $\{W_i \}_{i\in I}$ be a frame of subspaces with respect to $\{v_i\}_{i \in I}$
for $\H$. Then the frame operator $S_{W,v}$ for $\{W_i \}_{i\in I}$ and
$\{v_i\}_{i \in I}$ is defined by
\[S_{W,v}(f) = T_{W,v}T_{W,v}^*(f) = T_{W,v}(\{v_i{\pi}_{W_i}(f)\}_{i\in I}) =
\sum_{i\in I}v_i^2{\pi}_{W_i}(f).\]
\end{definition}

The next proposition generalizes a result of Fornasier \cite{F,F1}.

\begin{prop} \label{reconstruction}
Let $\{W_i \}_{i\in I}$ be a frame of subspaces with respect to $\{v_i\}_{i \in I}$
with frame bounds $C$ and $D$.
Then the frame operator $S_{W,v}$ for $\{W_i \}_{i\in I}$ and $\{v_i\}_{i \in I}$
is a positive, self-adjoint, invertible operator on $\H$ with $CI\le S_{W,v} \le
DI$. Further, we have the reconstruction formula
\[ f = T_{S_{W,v}^{-1}W,v} T_{W,v}^* (f) = \sum_{i \in I} v_i^2 S_{W,v}^{-1} \pi_{W_i}(f)
\quad \mbox{for all } f \in \H.\]
\end{prop}

\begin{proof}
For any $f\in  \H$, we have
\[\langle S_{W,v}(f), f \rangle  = \langle \sum_{i\in I}v_i^2{\pi}_{W_i}(f),
f\rangle = \sum_{i\in I}v_i^2\langle {\pi}_{W_i}(f),f \rangle = \sum_{i\in
I}v_i^2\|{\pi}_{W_i}(f)\|^2,\]
which implies that $S_{W,v}$ is a positive operator. We further compute
\[\langle Cf , f \rangle = C\|f\|^2 \le \sum_{i\in I}v_i^2\|{\pi}_{W_i}(f)\|^2
= \langle S_{W,v}(f), f\rangle\le \langle Df , f \rangle.\]
This shows that $CI\le S_{W,v}\le DI$ and hence $S_{W,v}$ is an invertible operator on
$\H$.
Furthermore, for any $f,g\in \H$ we have
\[\langle S_{W,v}(f), g\rangle = \sum_{i\in I} v_i^2\langle {\pi}_{W_i}(f), g
\rangle = \sum_{i\in I}v_i^2\langle f, {\pi}_{W_i}(g) \rangle.\]
Thus $S_{W,v}$ is self-adjoint.
At last the reconstruction formula follows immediately from
\[ f = S_{W,v}^{-1} S_{W,v}(f) = \sum_{i \in I} v_i^2 S_{W,v}^{-1} \pi_{W_i}(f).\]
\end{proof}

The following result will show the connection between the frame operator
for a frame of subspaces and the frame operator for the frame generated
by orthonormal bases of the subspaces. Also the connection between the
reconstruction formulas is exposed.

\begin{prop}\label{transfer_frameoperator}
Let $\{W_i \}_{i\in I}$ be a frame of subspaces with respect to $\{v_i\}_{i \in I}$ for $\H$ and
$\{v_i f_{ij} \}_{j\in J_i}$ be a Parseval frame for $W_i$ for each
$i\in I$.  Then the frame operator $S_{W,v}$
equals the frame operator $S_{vf}$ for the frame
$\{v_i f_{ij}\}_{i\in I,j\in J_i}$, and for all $g \in \H$ we have
\[ \sum_{i \in I} v_i^2 S_{W,v}^{-1} \pi_{W_i}(g)
= \sum_{i\in I, j\in J_i} \ip{g}{v_i f_{ij}} S_{vf}^{-1} v_i f_{ij}.\]
\end{prop}

\begin{proof}
Since $\{f_{ij}\}_{j\in J_i}$ is a Parseval frame for $W_i$ for all
$i\in I$, if $g\in \H$ then
\[{\pi}_{W_i}(g) = \sum_{j\in J_i}\langle {\pi}_{W_i}(g), f_{ij}\rangle f_{ij}
= \sum_{j\in J_i}\langle g, f_{ij}\rangle f_{ij}.\]
Thus
\[S_{W,v}(g) = \sum_{i\in I}v_i^2 {\pi}_{W_i}(g) = \sum_{i\in I, j\in J_i}
\langle g,v_i f_{ij}\rangle v_i f_{ij} = S_{vf}(g).\]
Moreover, we obtain
\[ \sum_{i \in I} v_i^2 S_{W,v}^{-1} \pi_{W_i}(g)
= \sum_{i \in I} S_{vf}^{-1} \sum_{j \in J_i} \ip{g}{v_i f_{ij}} v_i f_{ij}
= \sum_{i\in I, j\in J_i} \ip{g}{v_i f_{ij}} S_{vf}^{-1} v_i f_{ij}.\]
\end{proof}

Using the frame operator for a frame of subspaces for a special subspace
yields an easy way to compute the orthogonal projection onto this subspace.

\begin{prop}
Let $\{W_i\}_{i\in I}$ be a frame of subspaces with respect to $\{v_i\}_{i \in I}$
for a subspace $V$ of $\H$.
Then, the orthogonal projection $\pi_V$ onto $V$ is given by
\[ \pi_V(f) = \sum_{i \in I} v_i^2 S_{W,v}^{-1} \pi_{W_i}(f)\quad \mbox{for all } f \in \H.\]
\end{prop}

\begin{proof}
The fact that $S_{W,v} : V \to V$ implies that $\pi_V(f) = 0$ for all $f \in V^\perp$.
By Proposition \ref{reconstruction}, we have
\[ f =  \sum_{i \in I} v_i^2 S_{W,v}^{-1} \pi_{W_i}(f) \quad \mbox{for all } f \in V.\]
Thus $\pi_V^2=\pi_V$, which finishes the proof.
\end{proof}

In the same manner as in frame theory we define a dual frame of
subspaces associated with a frame of subspaces.

\begin{definition}
Let $\{W_i \}_{i\in I}$ be a frame of subspaces with respect to $\{v_i\}_{i \in I}$
and with frame operator $S_{W,v}$.
Then $\{S_{W,v}^{-1}W_i\}_{i\in I}$ is called the {\em dual frame of
subspaces with respect to $\{v_i\}_{i \in I}$}.
\end{definition}

The dual frame of subspaces is a frame of subspaces with the same weights.  In
fact, more is true.

\begin{prop}
Let $\{W_i\}_{i\in I}$ be a frame of subspaces with respect to $\{v_i\}_{i\in I}$,
and let $T:\H \rightarrow \H$ be an invertible operator on $\H$. Then $\{TW_i\}_{i\in I}$
is a frame of subspaces with respect to $\{v_i\}_{i\in I}$.
\end{prop}

\begin{proof}
Since $T$ is an invertible operator on $\H$, we have that
${\pi}_{TW_i} = T{\pi}_{W_i}T^{-1}$.  Let $C,D>0$ be the frame bounds for the
frame of subspaces $\{W_i\}_{i\in I}$.  Then for all $f\in \H$ we have
\[ \sum_{i\in I}v_i^2 \|{\pi}_{TW_i}(f)\|^2
= \sum_{i\in I}v_i^2 \|T{\pi}_{W_i}T^{-1}(f)\|^2
\le \|T\|^2 \sum_{i\in I}v_i^2 \|{\pi}_{W_i}T^{-1}(f)\|^2\]
\[ \le \|T\|^2 D \|T^{-1}(f)\|^2
\le \|T\|^2 \|T^{-1}\|^2 D\|f\|^2 . \]
Similarly, we obtain a lower frame of subspaces bound for $\{TW_i\}_{i\in I}$.
\end{proof}


\subsection{Parseval frames of subspaces}
\label{Parseval_fos}

Parseval frames play an important role in abstract frame theory, since they are
extremely useful for applications. Therefore in this subsection we study characterizations
of Parseval frames of subspaces and special cases of them.

The first result extends \cite[Corollary 4.1]{CK01}.

\begin{coro} \label{transfer_Parseval}
For each $i \in I$ let $v_i > 0$ and let $\{f_{ij}\}_{j \in J_i}$ be
a Parseval frame sequence in $\mathcal{H}$. Define
$W_i = \span_{j\in J_i}\{f_{ij}\}$ for all $i\in I$, and choose for
each subspace $W_i$ an orthonormal basis $\{e_{ij}\}_{j \in J_i}$.
Then the following conditions are equivalent.
\begin{enumerate}
\item $\{v_i f_{ij}\}_{i \in I,j \in J_i}$ is a Parseval frame for $\mathcal{H}$.
\item $\{v_i e_{ij}\}_{i \in I, j \in J_i}$ is a Parseval frame for $\H$.
\item $\{W_i\}_{i\in I}$ is a Parseval frame of subspaces with respect to $\{v_i\}_{i \in I}$
for $\H$.
\end{enumerate}
\end{coro}

\begin{proof}
This follows immediately from Theorem \ref{transfer_frame}.
\end{proof}

We can also characterize Parseval frames of subspaces in terms of their frame operators in
a similar manner as in frame theory.

\begin{prop} \label{prop_Parseval}
Let $\{W_i \}_{i\in I}$ be a family of subspaces in $\H$, and let $\{v_i\}_{i \in I}$
be a family of weights.  Then the following conditions are equivalent.
\begin{enumerate}
\item $\{W_i \}_{i\in I}$ is a Parseval frame of subspaces with respect to $\{v_i\}_{i \in I}$ for $\H$.
\item $S_{W,v} = I$.
\end{enumerate}
\end{prop}

\begin{proof}
For each $i \in I$, let $\{e_{ij}\}_{j \in J_i}$ be an orthonormal basis for $W_i$.

By Proposition \ref{reconstruction}, (1) implies (2). To prove the converse implication
suppose that $S_{W,v} = I$. Then for all $f \in \H$ we have
\[ f = S_{W,v}(f) = \sum_{i \in I} v_i^2 \pi_{W_i}(f) = \sum_{i \in I} v_i^2 \sum_{j \in J_i}
\ip{f}{e_{ij}} e_{ij}.\]
This yields
\[ \|f\|^2 = \ip{\sum_{i \in I} v_i^2 \sum_{j \in J_i} \ip{f}{e_{ij}} e_{ij}}{f}
= \sum_{i \in I} v_i^2 \|\pi_{W_i}(f)\|^2.\]
\end{proof}

We also have the following characterization of orthonormal bases of subspaces, which
reflects exactly the situation in frame theory.

\begin{prop} \label{prop_orthonormalbasis}
Let $\{W_i \}_{i\in I}$ be a family of subspaces in $\H$, and let $\{v_i\}_{i \in I}$
be a family of weights.  Then the following conditions are equivalent.
\begin{enumerate}
\item $\{W_i \}_{i\in I}$ is an orthonormal basis of subspaces for $\H$.
\item $\{W_i \}_{i\in I}$ is a $1$-uniform Parseval frame of subspaces for $\H$.
\end{enumerate}
\end{prop}

\begin{proof}
For each $i \in I$, let $\{e_{ij}\}_{j \in J_i}$ be an orthonormal basis for $W_i$.

If (1) is satisfied, then $\{e_{ij}\}_{i \in I, j \in J_i}$
is an orthonormal basis for $\H$. This implies
\[ \|f\|^2 =\sum_{i \in I} \sum_{j \in J_i}\absip{e_{ij}}{f}^2= \sum_{i \in I} \|\pi_{W_i}(f)\|^2\]
for all $f \in \H$. Thus also (2) holds.

On the other hand suppose that (2) holds. Then for all $f \in \H$ we have
\[ \|f\|^2 = \sum_{i \in I} \|\pi_{W_i}(f)\|^2 = \sum_{i \in I} \sum_{j \in J_i}
\absip{e_{ij}}{f}^2\]
and $\|e_{ij}\| = 1$ for all $i \in I,j \in J_i$, which shows that $\{e_{ij}\}_{i \in I, j \in J_i}$
is an orthonormal basis for $\H$. This immediately implies $\H = \bigoplus_{i \in I} W_i$, hence (1)
follows.
\end{proof}


\subsection{Resolution of the identity}
\label{roti}

Let $\{W_i \}_{i\in I}$ be a frame of subspaces with respect to $\{v_i\}_{i \in I}$
for $\H$ and let its frame operator be denoted by $S_{W,v}$. By Proposition
\ref{reconstruction}, we have
\[ f = \sum_{i \in I} v_i^2 S_{W,v}^{-1} \pi_{W_i}(f) \quad \mbox{for all } f \in \H.\]
This shows that the family of operators $\{v_i^2 S_{W,v}^{-1} \pi_{W_i}\}_{i \in I}$ is a
resolution of the identity. But a frame of subspaces for $\H$ provides us with
many more resolutions of the identity than only this one.

We start our consideration with the general definition of a resolution of the
identity.

\begin{definition}
Let $I$ be an indexing set.
A family of bounded operators $\{T_i\}_{i\in I}$ on $\H$
is called a {\em (unconditional) resolution of the identity}
on $\H$, if for all $f\in \H$ we have
\[f = \sum_{i\in I} T_{i}(f)\]
(and the series converges unconditionally for all $f\in \H$).
\end{definition}

Note that it follows from the definition and the uniform
boundedness principle that $\sup_{i\in I}\|T_{i}\|<\infty$.

The following result shows another way to obtain a resolution of the identity
from a frame of subspaces, which even satisfies an analog of (\ref{deffos}).

\begin{prop}
Let $\{v_i\}_{i \in I}$ be a family of weights,
and for each $i \in I$ let $\{v_i f_{ij}\}_{j \in J_i}$ be a frame sequence in $\mathcal{H}$
with frame bounds $A_i$ and $B_i$. Suppose that
$\{W_i\}_{i \in I}$ is a frame of subspaces with respect to $\{v_i\}_{i \in I}$
for $\H$ with frame bounds $C$ and $D$, where $W_i = \span_{j\in J_i}\{f_{ij}\}$
for all $i\in I$. Then $\{v_i f_{ij}\}_{i\in I,j\in J_i}$ is a frame for $\H$
with frame operator denoted by $S_{vf}$.
Further, for each $i \in I$, let $T_i :\H\rightarrow W_i$ be given by
\[T_i (f) = \sum_{j\in J_i}\langle f, S_{vf}^{-1}v_i f_{ij}
\rangle v_i f_{ij}.\]
If $0 < A = \inf_{i \in I} A_i \le B = \sup_{i \in I} B_i < \infty$,
then $\{T_i\}_{i \in I}$ is an unconditional resolution of the identity on $\H$
satisfying
\[\frac{AC}{B^2D^2}\|f\|^2 \le
\sum_{i\in I}v_i^2\|T_i(f)\|^2 \le \frac{B^2D^3}{A^2C^2}\|f\|^2 \quad
\mbox{for all } f\in \H.\]
\end{prop}

\begin{proof}
Recall that $\{v_i f_{ij}\}_{i\in I,j\in J_i}$ is a frame for $\H$ by Theorem
\ref{transfer_frame} with frame bounds $AC$ and $BD$.
For any $f\in \H$ we have
\[f =\sum_{i\in I} \sum_{j\in J_i}\langle f,S_{vf}^{-1}v_i f_{ij}\rangle v_i f_{ij}
= \sum_{i\in I}T_{i}(f).\]
Since this is convergence relative to a frame, the convergence is unconditional.

For each $i\in I$, let $S_{vf,i}$ be the frame operator for $\{v_i f_{ij}\}_{j\in
J_i}$. Let $i\in I$ be fixed. Then we obtain
\begin{eqnarray*}
\|T_i (f)\|^2 & = & \| \sum_{j\in J_i}\langle S_{vf}^{-1}f,v_i f_{ij}\rangle v_i f_{ij}\|^2\\
& = & \|S_{vf,i} \pi_{W_i}S_{vf}^{-1}(f)\|^2\\
& \le & \|S_{vf,i}\|^2 \|{\pi}_{W_i}S_{vf}^{-1}(f)\|^2.
\end{eqnarray*}
To prove the upper bound, we compute
\[\sum_{i\in I} v_i^2\|T_i (f)\|^2 \le B^2D^{2}\sum_{i\in
I}v_i^2 \|{\pi}_{W_i}S_{vf}^{-1}(f)\|^2
\le B^2D^3 \|S_{vf}^{-1}(f)\|^2 \le \frac{B^2D^3}{A^2C^2}\|f\|^2.\]
The lower bound follows from
\[\sum_{i\in I}v_i^2\|T_i(f)\|^2 = \sum_{i \in I}v_i^2 \|S_{vf,i} {\pi}_{W_i}S^{-1}(f)\|^2
\ge \sum_{i\in I}A_i v_i^2 \|{\pi}_{W_i}S_{vf}^{-1}(f)\|^2\]
\[\ge AC\|S_{vf}^{-1}(f)\|^2
\ge \frac{AC}{B^2D^2}\|f\|^2.\]
\end{proof}

We now give another method for obtaining an unconditional resolution of the identity
from a frame of subspaces.  A special case of this can be found in Fornasier \cite{F,F1}.

\begin{prop}
Let $\{W_i \}_{i\in I}$ be a frame of subspaces with respect to $\{v_i\}_{i \in I}$
for $\H$ with frame bounds $C$ and $D$, and let $S_{W,v}$ denote its frame operator.
Then $\{T_i\}_{i \in I}$ defined by $T_i = {\pi}_{W_i}S_{W,v}^{-1}$, $i \in I$ satisfies
that $\{v_i^2 T_i\}_{i \in I}$
is an unconditional resolution of the identity, and for all $f\in \H$ we have
\[\frac{C}{D^2}\|f\|^2
 \le \sum_{i\in I}v_i^2\|T_i (f)\|^2 \le \frac{D}{C^2}\|f\|^2.\]
\end{prop}

\begin{proof}
First, for any $f\in \H$ we have
\[\sum_{i\in I}v_i^2{\pi}_{W_i}S_{W,v}^{-1}(f) = S_{W,v}S_{W,v}^{-1}(f) = f.\]
To prove the second claim we compute
\[\frac{C}{D^2}\|f\|^2 \le
C\|S_{W,v}^{-1}(f)\|^2 \le \sum_{i\in I}v_i^2 \|{\pi}_{W_i}S_{W,v}^{-1}(f)\|^2
\le D\|S_{W,v}^{-1}(f)\|^2 \le \frac{D}{C^2}\|f\|^2.\]
\end{proof}

The next result will turn out to be useful for proving a lower bound for
$\sum_{i\in I}v_i^2\|T_{i}(f)\|^2$ if $\{v_i^2 T_i\}_{i \in I}$ is a resolution of the
identity.

\begin{lemma} \label{resolutionlower}
Let $\{W_i \}_{i\in I}$ be a frame of subspaces with respect to $\{v_i\}_{i \in I}$
with frame bounds $C$ and $D$ for $\H$, and let $T_i :\H
\rightarrow W_i$ be such that $\{v_i^2 T_i\}_{i \in I}$ is a resolution of
the identity on $\H$  (Note that a resolution
of the identity need not be unconditional so the index set must have an
ordering on it.  In our case, the result will hold for any ordering so
we do not specify the ordering here).  For any $J\subset I$ we have
\[\frac{1}{D}\|\sum_{j\in J}v_i^2 T_{j}(f)\|^2 \le \sum_{j\in J}v_j^2\|T_{j}(f)\|^2
\quad \mbox{for all } f \in \H.\]
\end{lemma}

\begin{proof}
We may assume that $|J| < \infty$, since if our inequality
holds for all finite subsets then it holds for all subsets.
Let $f \in \H$ and set $g = \sum_{j\in J}v_j^2 T_{j}(f)$. Then, using the
fact that $\{W_i \}_{i\in I}$ is a frame of subspaces with respect to $\{v_i\}_{i \in I}$
for $\H$, we compute
\begin{eqnarray*}
\|g\|^4 & = & \left(\langle g, \sum_{j\in J}v_j^2 T_{j}(f)\rangle \right)^2\\
& = & \left(\sum_{j\in J}v_j \langle g, v_j T_{j}(f) \rangle \right)^2\\
& = & \left( \sum_{j\in J}v_j\langle{\pi}_{W_j}(g), v_j T_{j}(f) \rangle \right)^2\\
& \le & \left( \sum_{j\in J}v_j\|{\pi}_{W_j}(g)\|\|v_j T_{j}(f)\| \right)^2\\
& \le & \sum_{j\in J}v_j^2 \|{\pi}_{W_j}(g)\|^2 \sum_{j\in J}\|v_j T_{j}(f)\|^2\\
& \le & D\|g\|^2 \sum_{j\in J}\|v_j T_{j}(f)\|^2.
\end{eqnarray*}
Dividing both sides of this inequality by $D\|g\|^2$ completes the proof.

\end{proof}

Using this lemma, we obtain bounds for $\sum_{i\in I}v_i^2\|T_{i}(f)\|^2$ for
many resolutions of the identity $\{v_i^2 T_i\}_{i \in I}$.

\begin{prop} \label{Tpi=T}
Let $\{W_{i}\}_{i\in I}$ be a frame of subspaces with respect to $\{v_i\}_{i \in I}$ for
$\H$ with frame bounds $C$ and $D$, let $T_{i}:\H \rightarrow W_i$ be such that
$\{v_i^2 T_i\}_{i \in I}$ is a resolution of the identity on $\H$,
and assume that $T_{i}{\pi}_{W_{i}} = T_i$.  Then for all $f\in \H$
\[\frac{1}{D}\|f\|^2 \le \sum_{i\in I}v_i^2\|T_i (f)\|^2 \le DE\|f\|^2,\]
where $E= \mbox{sup}_{i}\|T_i \| < \infty$.
\end{prop}

\begin{proof}
By Lemma \ref{resolutionlower}, for all $f \in \H$, we have
\[\frac{1}{D}\|f\|^2 =\frac{1}{D}\| \sum_{i\in I}v_i^2 T_{i}(f)\|^2 \le \sum_{i\in I}v_i^2\|T_{i}(f)\|^2 =
\sum_{i\in I}v_i^2\|T_{i}{\pi}_{W_i}(f)\|^2\]
\[\le \sum_{i\in I}v_i^2\|T_{i}\|^2 \|{\pi}_{W_i}(f)\|^2
 \le E\sum_{i\in I}v_i^2\|{\pi}_{W_i}(f)\|^2 \le DE\|f\|^2.\]
\end{proof}

Obviously the condition $T_{i}{\pi}_{W_{i}} = T_i$ for all $i \in I$ is satisfied
by the example $\{v_i^2 S_{W,v}^{-1} \pi_{W_i}\}_{i \in I}$ from the beginning of this
subsection. This shows that this family of operators is not only a resolution of
the identity but even satisfies an analog of (\ref{deffos}).

The following definition provides us with a condition, which
implies that a resolution of the identity $T_{i}:\H \rightarrow W_i$, where
$\{W_{i}\}_{i\in I}$ is a frame of subspaces with respect to $\{v_i\}_{i \in I}$
for $\H$, automatically satisfies an analog of (\ref{deffos}).

\begin{definition}
A family of bounded operators $\{T_{i}\}_{i\in I}$ on
$\H$ is called an {\em ${\ell}^{2}$-resolution of the identity with
respect to a family of weights $\{v_i\}_{i \in I}$} on $\H$, if it is a resolution of
the identity on $\H$ and there exists a constant $B>0$ so that for
all $f\in \H$ we have
\[\sum_{i\in I}v_i^{-2}\|T_{i}(f)\|^2 \le B\|f\|^2.\]
\end{definition}

\begin{theorem}
Let $\{W_{i}\}_{i\in I}$ be a frame of subspaces with respect to $\{v_i\}_{i \in I}$ for $\H$,
and let $T_{i}:\H \rightarrow W_i$ be such that $\{v_i^2 T_i\}_{i \in I}$ is an
${\ell}^2$-resolution of the identity with respect to $\{v_i\}_{i \in I}$ on $\H$.
Then there exist constants $A,B>0$ so that for all $f\in \H$ we have
\[A\|f\|^2 \le \sum_{i\in I}v_i^2\|T_{i}(f)\|^2 \le B\|f\|^2.\]
\end{theorem}

\begin{proof}
This follows immediately from the definition of an ${\ell}^2$-resolution of the
identity on $\H$ and Lemma \ref{resolutionlower}.
\end{proof}


\section{Riesz decompositions}
\label{Riesz_dec}

In this section we first study minimal frames of subspaces, which share
similar properties with minimal frames.

\begin{definition}
A family of subspaces $\{W_i\}_{i\in I}$ of $\H$
is called {\em minimal}, if for each $i \in I$
\[W_i \cap \span_{j \in I, j \neq i} \{W_j\} = \{0\}.\]
\end{definition}

Using orthonormal bases for the subspaces we obtain a useful characterization of
minimal families of subspaces.

\begin{lemma} \label{transfer_minimal}
Let $\{W_i\}_{i\in I}$ be a family of subspaces in $\H$, and for each $i \in I$
let $\{e_{ij}\}_{j \in J_i}$ be an orthonormal basis for $W_i$.
Then the following conditions are equivalent.
\begin{enumerate}
\item $\{W_i\}_{i\in I}$ is minimal.
\item $\{e_{ij}\}_{i \in I, j \in J_i}$ is minimal.
\end{enumerate}
\end{lemma}

\begin{proof}
The implication (2) $\Rightarrow$ (1) is obvious.
To prove (1) $\Rightarrow$ (2) suppose that $\{c_{ij}\}_{j\in J_i}\in {\ell}_{2}(J_i )$
for all $i\in I$ and we have
$f_i = \sum_{j\in J_i}c_{ij}e_{ij}$ and $f = \{f_i \}_{i\in I}
\in \left ( \sum_{i\in I}\oplus W_i \right ) _{{\ell}_2}$.
If $\sum_{i\in I}f_i = 0$, then by minimality of $\{W_i \}_{i\in I}$ we have
that $f_i = 0$ for all $i\in I$ and so $c_{ij}=0$ for all $i\in I, j\in J_i$.  It
follows that $\{e_{ij}\}_{i\in I, j\in J_i}$ is a minimal frame for
$\H$.
\end{proof}

The following two propositions show that for families of subspaces we can also give a
definition of biorthogonal families of subspaces, which possess similar properties compared
to minimal frames of subspaces as in the situation of minimal frames
(compare \cite[Lemma 3.3.1]{Chr03}).

\begin{prop}
Let $\{W_i\}_{i\in I}$ be a family of subspaces in $\H$.
Then the following conditions are equivalent.
\begin{enumerate}
\item $\{W_i\}_{i\in I}$ is minimal.
\item There exists a unique maximal (up to containment)
biorthogonal family of subspaces for $\{W_i\}_{i\in I}$,
i.e., there exists a family of subspaces $\{V_i\}_{i\in I}$ with
$W_i \perp V_j$ for all $i,j \in I$, $j \neq i$ and $f \not\perp V_i$ for all
$f \in W_i$, $i \in I$.
\end{enumerate}
Moreover, if $\{W_i\}_{i\in I}$ is a minimal frame of subspaces with respect to
$\{v_i\}_{i \in I}$ for $\H$, then $\{S_{W,v}^{-1/2}W_i \}_{i\in I}$ is an orthogonal family
of subspaces in $\H$.
\end{prop}

\begin{proof}
Suppose that (2) holds and towards a contradiction assume that there exists
$i \in I$ and $0 \neq f = \sum_{j \in I, j \neq i} g_j \in W_i$ with
$g_j \in W_j$. By (2), we have $g_j \perp V_i$ for all $j \neq i$,
hence $f \perp V_i$, but this is a contradiction. Thus (1) follows.

To prove the opposite direction suppose that $\{W_i\}_{i\in I}$ is minimal.
For each $i \in I$, let $P_i$ denote the
orthogonal projection onto $\span_{j \in I, j \neq i} \{W_j\}$. Let $i \in I$ and let $V_i$
be defined by $V_i = (I-P_i) \H$ for all $i\in I$.
By the definition of $V_i$, we have $W_j \perp V_i$ for all $j \neq i$. Towards
a contradiction assume that
there exists $f \in W_i$ with $\ip{f}{g} = 0$ for all $g \in V_i$. Then
$f\in P_{i} \H$ and so $W_{i}\cap P_{i}\H \not= \{0\}$, which is a contradiction.

For the moreover part, let $\{v_i e_{ij}\}_{j\in J_i}$ be an orthonormal basis for
$W_i$ for each $i\in I$.  By Proposition \ref{transfer_frameoperator}, $S_{W,v}$ is the frame operator
for $\{v_i e_{ij}\}_{i\in I, j\in J_i}$. Since $\{W_i\}_{i\in I}$ is minimal,
Lemma \ref{transfer_minimal} implies that $\{v_i e_{ij}\}_{i\in I, j\in J_i}$ is
a minimal frame for $\H$ and hence is a Riesz basis for $\H$.  Thus
$\{S_{W,v}^{-1/2}v_i e_{ij}\}_{i\in I,j\in J_i}$ is an orthonormal sequence in $\H$.
Since we have
\[S_{W,v}^{-1/2}W_i = \span_{j\in J_i} \{S_{W,v}^{-1/2}v_i e_{ij}\},\]
it follows that $\{S_{W,v}^{-1/2}W_i \}_{i\in I, j\in J_i }$ is an orthogonal sequence
in $\H$.
\end{proof}

The next definition transfers the definition of Riesz bases and exact sequences
to families of subspaces in a canonical way. The so-called Riesz decomposition
will share most properties with Riesz bases. However it will turn
out that being an exact frame of subspaces is strictly weaker than being a
Riesz decomposition, a fact which differs from the situation in abstract
frame theory.

\begin{definition}
We call a frame  of subspaces $\{W_i \}_{i\in I}$ with respect to some family of weights
for $\H$ a {\em Riesz decomposition} of $\H$, if
for every $f\in \H$ there is a unique choice of $f_i \in W_i $ so that $f=\sum_{i\in I} f_i$.
A frame of subspaces with respect to some family of weights is {\em exact}, if it ceases
to be a frame of subspaces once one element is deleted.
\end{definition}

\begin{lemma} \label{transfer_Riesz_exact}
Let $\{W_i\}_{i\in I}$ be a family of subspaces in $\H$, and for each $i \in I$
let $\{e_{ij}\}_{j \in J_i}$ be an orthonormal basis for $W_i$.
\begin{enumerate}
\item The following conditions are equivalent.
\begin{enumerate}
\item $\{W_i\}_{i\in I}$ is a Riesz decomposition of $\H$.
\item $\{e_{ij}\}_{i \in I, j \in J_i}$ is a Riesz basis for $\H$.
\item $\{e_{ij}\}_{i \in I, j \in J_i}$ is an unconditional basis for
$\H$.
\end{enumerate}
\item Let $\{W_i\}_{i\in I}$ be a frame of subspaces with respect to $\{v_i\}_{i \in I}$ for $\H$.
If $\{v_i e_{ij}\}_{i \in I, j \in J_i}$ is an exact frame for $\H$, then also $\{W_i\}_{i\in I}$
is an exact frame of subspaces with respect to $\{v_i\}_{i \in I}$ for $\H$. The opposite implication
is not valid.
\end{enumerate}
\end{lemma}

\begin{proof}
First we prove (1). The equivalence (b) $\Leftrightarrow$ (c) follows immediately from the
fact that $\{v_i e_{ij}\}_{i \in I, j \in J_i}$
is bounded and that a Schauder basis is a Riesz basis if and only if it is a bounded unconditional
basis. The implication (b) $\Rightarrow$ (a) is obvious.
To prove (a) $\Rightarrow$ (c) assume that $\{e_{ij}\}_{i \in I, j \in J_i}$ is not an
unconditional basis. Hence there exist $f \in \H$ and
sequences $\{c_{ij}\}_{i \in I, j \in J_i}$ and $\{d_{ij}\}_{i \in I, j \in J_i}$
with $f = \sum_{i \in I, j \in J_i} c_{ij} e_{ij} = \sum_{i \in I, j \in J_i} d_{ij} e_{ij}$
such that $c_{i_0j_0} \neq d_{i_0j_0}$ for some $i_0 \in I$, $j_0 \in J_{i_0}$. By construction
$\{e_{i_0j}\}_{j \in J_{i_0}}$ is an orthonormal basis for $W_{i_0}$, hence
$\sum_{j \in J_{i_0}} c_{i_0j} e_{i_0j} \neq \sum_{j \in J_{i_0}} d_{i_0j} e_{i_0j}$,
which implies that $\{W_i\}_{i\in I}$ is not a Riesz decomposition.

To prove (2) suppose that $\{W_i\}_{i\in I}$ is a frame of subspaces with respect to
$\{v_i\}_{i \in I}$ for $\H$. By Theorem \ref{transfer_frame}, $\{v_i e_{ij}\}_{i \in I, j \in J_i}$
is a frame for $\H$. If this is exact, then, by definition, deleting one element $v_{i_0} e_{i_0j_0}$
does not leave a frame. Thus also $\{v_i e_{ij}\}_{i \in I \backslash \{i_0\}, j \in J_i}$
does not form a frame. Applying  Theorem \ref{transfer_frame} once more yields the first claim.
The fact that the opposite implication is not valid
is demonstrated by Example \ref{exa_exact}.
\end{proof}

The next theorem is the analog of a well-known result in abstract frame
theory (see \cite[Theorem 6.1.1]{Chr03}), only the role of exact frames
of subspaces differs from the frame situation.

\begin{theorem} \label{Riesz_minimal_exact}
Let $\{W_i \}_{i\in I}$ be a frame of subspaces with respect to $\{v_i\}_{i \in I}$
for $\H$. Then the following conditions are equivalent.
\begin{enumerate}
\item $\{W_i \}_{i\in I}$ is a Riesz decomposition of $\H$.
\item $\{W_i \}_{i\in I}$ is minimal.
\item The synthesis operator $T_{W,v}$ is one-to-one.
\item The analysis operator $T_{W,v}^*$ is onto.
\end{enumerate}
Moreover, if $\{W_i \}_{i\in I}$ is a Riesz decomposition of $\H$, then
it is also an exact frame of subspaces for $\calH$. The opposite implication
is not valid.
\end{theorem}

\begin{proof}
First note that (3) $\Leftrightarrow$ (4) is always true for operators on
a Hilbert space. Moreover, it is obvious that (1) implies (3).

Next we prove (4) $\Rightarrow$ (1).  By Theorem \ref{PP1}, $T_{W,v}^*$ is
an isomorphism.  Therefore if it is onto, then it is invertible.
Hence, $T_{W,v}$ is invertible.  This implies that for every $f\in \H$ there exists
a $\{f_i \}_{i\in I} \in \left(\sum_{i\in I}\oplus W_i\right)_{l^2}$ so that
\[f= T_{W,v}(\{f_i \}_{i\in I}) = \sum_{i\in I}v_if_i.\]
If we have $f=\sum_{i \in I}f_i = \sum_{i \in I}g_i$ with $ \{f_i\},\{g_i \}\in
\left(\sum_{i \in I} \oplus W_i\right)_{l^2}$,
then it follows that $T_{W,v}( \{v_i^{-1} f_i \}_{i \in I}) = T_{W,v}(\{v_i^{-1} g_i\}_{i \in I})$.
Since $T_{W,v}$ is one-to-one, $\{v_i^{-1} f_i \}_{i \in I} =
\{v_i^{-1} g_i\}_{i \in I}$ and so $f_i = g_i$ for all $i\in I$. This shows
the equivalence of (1), (3), and (4).

It remains to prove that (1) is equivalent to (2).
If $\{W_i \}_{i\in I}$ is not a Riesz decomposition of $\H$, there
exists an element $f \in \H$ and $f_i,g_i \in W_i$, $i \in I$ with $f_{i_0}
\neq g_{i_0}$ for some $i_0 \in I$ and $f = \sum_{i \in I} f_i = \sum_{i \in I} g_i$.
It follows $0 \neq g_{i_0} - f_{i_0} = \sum_{i \in I, i \neq i_0} (f_i - g_i)$
and $g_{i_0} - f_{i_0} \in W_{i_0}$. This proves
\[g_{i_0} - f_{i_0} \in W_{i_0} \cap \span_{i \in I, i \neq i_0}\{W_i\},\]
which implies that $\{W_i \}_{i\in I}$ is not minimal.

To prove the converse implication assume that $\{W_i \}_{i\in I}$
is not minimal. Then, for some $i_0 \in I$, there exists
$0 \neq f = \sum_{i \in I, i \neq i_0} f_i \in W_{i_0}$ with $f_i \in W_i$. Hence
\[0 = \sum_{i \in I, i \neq i_0} (f_i - f) = \sum_{i \in I} 0.\]
Thus $\{W_i \}_{i\in I}$  is not a Riesz decomposition of $\H$.

To prove the moreover part, suppose $\{W_i \}_{i\in I}$ is a frame of subspaces
with respect to $\{v_i\}_{i \in I}$ for $\H$ and a Riesz decomposition of $\H$.
We will prove that this implies that $\{W_i \}_{i\in I}$ is exact.
For this, fix some $i_0 \in I$. Since $T_{W,v}^*$ is onto, there exists an element
$f \in \calH$ such that $\pi_{W_{i_0}}(f) \neq 0$, but $\pi_{W_{i}}(f)=0$
for all $i \neq i_0$. Thus
\[ \sum_{i\in I} v_i^2 \|{\pi}_{W_i}(f)\|^2 = v_{i_0}^2 \|{\pi}_{W_{i_0}}(f)\|^2.\]
Therefore it is not possible to delete $W_{i_0}$ from the frame of
subspaces yet leave a frame of subspaces. Since $i_0$ was chosen
arbitrarily, the claim follows.

The fact that the opposite implication is not valid is demonstrated by Example
\ref{exa_exact}.
\end{proof}

Note that we could have also proven the equivalence of (1) and (2) and the
moreover part of the previous result by using Lemma \ref{transfer_minimal},
Lemma \ref{transfer_Riesz_exact}, and \cite[Theorem 6.1.1]{Chr03}.
We have chosen to add the extended proof here in order to enlighten the
use of frames of subspaces.

Next we give an example for the different role exactness
plays in the situation of families of subspaces.

\begin{example} \label{exa_exact}
Let $\{e_i\}_{i \in \ZZ}$ be an orthonormal basis for some Hilbert space $\H$
and define the subspaces $W_1, W_2$ by
\[ W_1 = \span_{i \ge 0} \{e_i\} \quad \mbox{and} \quad W_2 = \span_{i \le 0} \{e_i\}.\]
Then $\{W_1,W_2\}$ is a frame of subspaces with respect to weights $\{v_1,v_2\}$
with $v_1=v_2=:v > 0$, since
\[ v_1\|\pi_{W_1}(f)\|^2 + v_2\|\pi_{W_2}(f)\|^2
= \sum_{i \in \ZZ} v\absip{f}{e_{i}}^2 + v\absip{f}{e_0}^2
= v\|f\|^2 + v\absip{f}{e_0}^2\]
and
\[ v\|f\|^2 \le v\|f\|^2 + v\absip{f}{e_0}^2 \le 2 v\|f\|^2.\]
It is also exact, since when we delete one subspace the remaining one does not
span the space. But it is not a Riesz decomposition, because we can write the
element $e_0$ as $e_0 = e_0 + 0$ and $e_0 = 0 + e_0$. Thus the decomposition is
not unique.
Also observe that the sequence $\{v e_i\}_{i \ge 0} \cup \{v e_i\}_{i \le 0}$
is a frame, but is not exact.
\end{example}

We conclude this subsection by mentioning that orthonormal bases of subspaces
are special cases of Riesz decompositions.

\begin{coro}
If $\{W_i \}_{i\in I}$ is an orthonormal basis of subspaces for $\calH$, then it
is also a Riesz decomposition of $\calH$.
\end{coro}

\begin{proof}
This follows immediately from the definition of a Riesz decomposition and
Proposition \ref{prop_orthonormalbasis}.
\end{proof}


\section{Several constructions}
\label{constr}

In this section we will discuss several constructions concerning frames of subspaces,
frames, and Riesz frames. Recall that in addition to what follows we are already
equipped with some constructions by Theorem \ref{transfer_frame}, Corollary
\ref{transfer_Parseval}, and Lemma \ref{transfer_Riesz_exact}.

\subsection{Constructions of frames of subspaces}
\label{constr_fos}

Dealing with Bessel families of subspaces is important, since there are
easy ways to turn such a family into a frame of subspaces. One way is to
just add the subspace $W_{0} = \H$ to the family. Another more careful
method is the following one: Take any orthonormal basis for $\H$ and
partition its elements into the subspaces $W_i$, $i \in I$. Then add
the subspaces spanned by the remaining elements to the Bessel family.
This yields a frame of subspaces.

Using the synthesis operator $T_{W,v}$, we obtain a characterization
of Bessel sequences of subspaces.

\begin{prop} \label{Bessel}
Let $\{W_i \}_{i\in I}$ be a family of subspaces of $\H$, and let $\{v_i\}_{i \in I}$
be a family of weights. Then the following conditions are equivalent.
\begin{enumerate}
\item $\{W_i\}_{i \in I}$ is a Bessel sequence of subspaces with respect to $\{v_i\}_{i \in I}$
for $\H$.
\item The synthesis operator $T_{W,v}$ is bounded and linear.
\end{enumerate}
\end{prop}

\begin{proof}
First suppose that (1) holds. Then Lemma \ref{convergence}
shows that the series in the definition of the synthesis operator $T_{W,v}$
converges unconditionally. Moreover, we have
\[\sum_{i\in I}v_i^2\|{\pi}_{W_i}(f)\|^2 \le B\|f\|^2.\]
By definition of the analysis operator $T_{W,v}^*$,
\begin{equation}\label{def_T}
\|T_{W,v}^*(f)\|^2 = \sum_{i\in I}v_i^2\|{\pi}_{W_i}(f)\|^2.
\end{equation}
Since $\{W_i\}_{i \in I}$ is a Bessel sequence of subspaces with respect to
$\{v_i\}_{i \in I}$, this implies that $T_{W,v}^*$ is bounded. Hence also $T_{W,v}$
is bounded, which shows (2).

If (2) holds, then also $T_{W,v}^*$ is a bounded operator. This fact together with
(\ref{def_T}) yields (1).
\end{proof}

One possible application for frames of subspaces is to the problem of classifying
those $g\in L^{2}({\mathbb R})$ and $0< a,b\le 1$ so that $(g,a,b)$ yields
a Gabor frame (see Example \ref{exa_Gabor} below).  This is an exceptionally deep problem
even in the case of characteristic functions \cite{CK,J}.  But we have simple classifications
of when $\{e^{2{\pi}imbt}g(t)\}_{m\in \ZZ}$ is a frame sequence in $L^{2}({\mathbb R})$
and when $\{e^{2{\pi}imbt}g(t-na)\}_{m,n\in \ZZ}$ has dense span in $L^{2}({\mathbb R})$.
By our results, this family will be a Gabor frame for $L^{2}({\mathbb R})$ if and only
if $\{W_n \}_{n\in \ZZ}$ is a frame of subspaces where $W_n$ is the closed linear span of
$\{e^{2{\pi}imbt}g(t-na)\}_{m\in \ZZ}$.

For some applications, we would like to take
a frame for $\H$ and divide it into subsets so that the closed linear span of these
subsets is a frame of subspaces for $\H$.  This is not always possible.
But the next proposition shows that one of the needed inequalities
will always hold.

\begin{prop}\label{P17}
Let $\{f_j \}_{j\in J}$ be a frame for $\H$ with
frame bounds $A$ and $B$.  Let $\{J_i\}_{i\in I}$ be a partition of the indexing set
$J$, and for all $i \in I$ let $W_i$ denote the closed linear span of
$\{f_j \}_{j\in J_i}$.  Then for all $f\in \H$ we have
\[\frac{A}{B}\|f\|^2 \le \sum_{i\in I}\|{\pi}_{W_i}(f)\|^2.\]
Hence, if $|I|<\infty$, then $\{W_i\}_{i\in I}$ is a $1$-uniform frame
of subspaces for $\H$.
\end{prop}

\begin{proof}
We compute
\[A\|f\|^2 \le \sum_{j\in J}|\langle f, f_j \rangle |^2 =
\sum_{i\in I}\sum_{j\in J_i}|\langle f   , f_j \rangle |^2
= \sum_{i\in I}\sum_{j\in J_i}|\langle {\pi}_{W_i}(f), f_j \rangle |^2.\]
Recall that if a family of vectors is a $B$-Bessel family then every subfamily
is also $B$-Bessel. Thus
\[ \sum_{i\in I}\sum_{j\in J_i}|\langle {\pi}_{W_i}(f), f_j \rangle |^2
\le \sum_{i\in I}B\|{\pi}_{W_i}(f)\|^2.\]
This proves the first claim.

If $|I|<\infty$, then we have
\[\sum_{i\in I}\|{\pi}_{W_i}(f)\|^2 \le |I| \cdot \|f\|^2.\]
Hence in this case $\{W_{i}\}_{i\in I}$ is always a frame of subspaces for $\H$,
in particular a $1$-uniform frame of subspaces.
\end{proof}

An easy way to obtain a frame of subspaces is provided by the next result.

\begin{prop} \label{finitepartition}
Let $\{f_j\}_{j \in J}$ be a frame for $\H$, let $J = J_1 \cup \ldots
\cup J_n$ be a finite partition of $J$, and let $\{v_i\}_{i=1}^n$ be a family
of weights. Then $\{W_i\}_{i=1}^n$ is a frame
of subspaces with respect to $\{v_i\}_{i=1}^n$ for $\H$, where
$W_i = \span_{j \in J_i} \{f_j\}$.
\end{prop}

\begin{proof}
Let $f \in \H$. Obviously, $\|\pi_{W_i}(f) \|^2 \le \|f\|^2$ for all $1 \le i \le n$,
which implies that
\[ \sum_{i=1}^n v_i^2 \|\pi_{W_i}(f) \|^2 \le \max_{i=1,\ldots,n}\{v_i^2\} \cdot \|f\|^2.\]
Thus  $\{W_i\}_{i=1}^n$ is a Bessel sequence of subspaces with respect to
$\{v_i\}_{i=1}^n$. The lower bound follows from an application of Proposition \ref{P17}.
That is, for any $f\in \H$ we have
\[ \frac{A}{B}\|f\|^2
\le \sum_{i\in I}\|{\pi}_{W_i}(f)\|^2
\le \frac{1}{\max_{i=1, \ldots,n}\{v_i^2\}}\sum_{i\in I}v_i^2 \|{\pi}_{W_i}(f)\|^2.\]
\end{proof}

This partition of the frame elements is not always a partition into frame
sequences. Let us consider the case of Gabor systems. In the following example
we will show that a large class of Gabor systems can be written as a frame of
subspaces. Moreover, we can characterize those Gabor atoms, for which this
partition is a partition into frame sequences.

\begin{example} \label{exa_Gabor}
For each $a \in \mathbb{R}$, let the unitary operators $E_a, T_a$ on $L^2(\mathbb{R})$
be defined by
\[ E_a f(x) = e^{2 \pi i a x} f(x)\quad \mbox{and}\quad T_a f(x) = f(x-a).\]
Given a function $g \in L^2(\mathbb{R})$ and $a,b > 0$, the Gabor system
determined by $g$ and $a,b$ is defined by
\[ \calG(g,a,b) = \{E_{ma} T_{nb} g : m,n \in \ZZ\}.\]
Let
\[ Z : L^2(\RR) \to L^2([0,1)^2), \quad Zf(x,y) = \sum_{k \in \ZZ} f(x+k)e^{2 \pi i k y}.\]
denote the Zak transform on $L^2(\RR)$ (compare \cite{Jan88}).

Let $h \in L^2(\mathbb{R})$ and $q \in \mathbb{N}$. In the following we will consider
some Gabor system $\calG(h,a,b)$ with $a,b > 0$, $ab = \frac{1}{q}$. Using a
metaplectic transform it can be shown that this system is unitarily equivalent to
$\calG(g,\frac{1}{q},1)$ for some $g \in L^2(\RR)$ \cite[Proposition 9.4.4]{Gro01},
hence it suffices to consider this system. Now the
Gabor system $\calG(g,\frac{1}{q},1)$ in turn can be decomposed using the partition
\begin{equation} \label{Gabor}
\calG(g,\tfrac{1}{q},1) =
\bigcup_{j=0}^{q-1} \{E_{\frac{1}{q}(mq+j)}T_ng\}_{m,n \in \mathbb{Z}}.
\end{equation}
By Proposition \ref{finitepartition}, the set of the subspaces
$W_j := \span_{m,n \in \mathbb{Z}}\{E_{\frac{1}{q}(mq+j)}T_ng\}$, $j = 0,\ldots,q-1$
is indeed a frame of subspaces.

\medskip

In a second step we will investigate, whether the sequences
$\{E_{\frac{1}{q}(mq+j)}T_ng\}_{m,n \in \mathbb{Z}}$ are frame sequences. It will turn
out that this will not happen unless the Zak transform is discontinuous.

We first observe that the following conditions are equivalent.
\begin{enumerate}
\item The sequence $\{E_{\frac{1}{q}(mq+j)}T_ng\}_{m,n \in \mathbb{Z}}$ is a frame sequence
for each $0 \le j \le q-1$.
\item There exist $0 < A \le B < \infty$ such that
\[A \le |Zg(x,y)|^2 \le B\quad \mbox{for almost all }(x,y) \in [0,1)^2 \backslash V,\]
where $V = \{(x,y) \in [0,1)^2 : Zg(x,y) = 0\}$.
\end{enumerate}

\medskip

This can be proven as follows.
First notice that since
\[E_{\frac{1}{q}}(\{E_{\frac{1}{q}(mq+j)}T_ng\}_{m,n \in \mathbb{Z}})
=  \{E_{\frac{1}{q}(mq+j+1)}T_ng\}_{m,n \in \mathbb{Z}}\quad \mbox{for all } 0 \le j < q-1\]
and
\[E_{\frac{1}{q}}(\{E_{\frac{1}{q}(mq+q-1)}T_ng\}_{m,n \in \mathbb{Z}})
=  \{E_{\frac{1}{q}(mq)}T_ng\}_{m,n \in \mathbb{Z}}
= \{E_{m}T_ng\}_{m,n \in \mathbb{Z}},\]
the fact that $E_{\frac{1}{q}}$ is a unitary operator implies that condition (1) holds if and only if
$\{E_{m}T_ng\}_{m,n \in \mathbb{Z}}$ is a frame sequence.
Since $Z : L^2(\mathbb{R}) \to L^2([0,1)^2)$ is an isomorphism \cite{Jan88}
and it is an easy calculation to show that
$Z (E_mT_ng)(x,y) = E_m(x) E_n(y) Zg(x,y)$, condition (1) holds if and only if
$\{E_m(x) E_n(y) Zg(x,y)\}_{m,n \in \mathbb{Z}}$ is a frame sequence.
Since $\{E_m E_n\}_{m,n \in \mathbb{Z}}$
is an orthonormal basis for $L^2([0,1)^2)$, for each $f \in L^2([0,1)^2)$ we obtain
\[ \sum_{m,n \in \mathbb{Z}} \absip{f}{E_m E_n Zg}^2
= \norm{f \cdot \overline{Zg}}^2.\]
This implies that (1) is equivalent to
\[ A \norm{f}^2 \le \norm{f \cdot \overline{Zg}}^2 \le B \norm{f}^2 \quad \mbox{for all }
f \in \span_{m,n \in \mathbb{Z}} \{E_m(x) E_n(y) Zg(x,y)\},\]
which holds if and only if
\[ A \norm{f \cdot Zg}^2 \le \norm{f \cdot |Zg|^2}^2 \le B \norm{f \cdot Zg}^2 \quad \mbox{for all }
f \in L^2([0,1)^2).\]
It is easy to check that this is equivalent to (2), which proves the claim.

\medskip

Finally,  we consider $g \in L^2(\RR)$ with $Zg$ being continuous.
By \cite{Jan88}, this implies that $Zg$ has a zero. Hence condition (2) can never
be fulfilled. This shows that the sequences
$\{E_{\frac{1}{q}(mq+j)}T_ng\}_{m,n \in \mathbb{Z}}$ can only be frame sequences,
if the Zak transform $Zg$ is discontinuous.
\end{example}


\subsection{Constructions of frames and Riesz frames}
\label{constr_frame_Riesz}

If we have a frame for $\H$, using a frame of subspaces for $\H$ we can
construct new frames for $\H$ from these components.

\begin{prop}
Let $\{W_i \}_{i\in I}$ be a frame of subspaces with respect to $\{v_i\}_{i \in I}$
for $\H$ and let $\{f_j \}_{j\in J}$ be a frame for $\H$.  Then there exist $A,B>0$
so that $\{{\pi}_{W_i}S_{W,v}^{-1}(f_j) \}_{j\in J}$ is a frame
for $W_i$ with frame bounds $A$ and $B$ for each $i \in I$.  Hence
$\{{\pi}_{W_i}S_{W,v}^{-1}(f_j) \}_{i\in I,j\in J}$ is also a frame for $\H$.
\end{prop}

\begin{proof}
Since $S_{W,v}^{-1}$ is an invertible operator on $\H$ and $\{f_j \}_{j\in J}$
is a frame for $\H$, we have that $\{S_{W,v}^{-1}f_j \}_{j\in J}$ is a frame for
$\H$ with frame bounds $A$ and $B$.  Therefore $\{{\pi}_{W_i}S_{W,v}^{-1}(f_j) \}_{j\in J}$
is a frame for $W_j$ with frame bounds $A$ and $B$ for every $i\in I$.
Since $\{W_i \}_{i\in I}$ is a frame of subspaces for $\H$, we have that
$\{{\pi}_{W_i}S_{W,v}^{-1}(f_j) \}_{i\in I,j\in J}$ is a frame for $\H$ by
Theorem \ref{transfer_frame}.
\end{proof}

To construct Riesz frames for $\H$ we first need to give an analog
definition for families of subspaces.

\begin{definition}
We call a frame of subspaces $\{W_i\}_{i\in I}$ a {\em Riesz frame of
subspaces with respect to $\{v_i\}_{i \in I}$,} if there exist constants $C, D>0$ so that every subfamily
$\{W_i\}_{i\in J}$ with $J\subset I$ is a frame of subspaces with respect to $\{v_i\}_{i \in J}$
for its closed linear span with frame bounds $C$ and $D$.
\end{definition}

First we may ask whether subfamilies of a frame of subspaces are
automatically frames of subspaces for their closed linear spans. The
following example shows that this is not always the case.

\begin{example}
In general, if $\{W_i\}_{i\in I}$ is a $1$-uniform frame of subspaces and
$J\subset I$, then $\{W_i\}_{i\in J}$ need not be a frame of subspaces
for its closed linear span.  For example, let $\{e_i\}_{i=1}^{\infty}$
be an orthonormal basis for $\H$ and for each $i \in I$ define the subspaces $W_i^1$,
$W_i^2$, and $W_i^3$ by
\[W_i^1 = \ospan \{e_{2i}+\tfrac{1}{i}e_{2i+1}\}, \quad
W_i^2 = \ospan \{e_{2i}\},\quad \mbox{and} \quad
W_i^3 = \ospan \{e_{2i+1}\}.\]
Then it is easily checked that $\{W_i^1, W_i^2, W_i^3 \}_{i=1}^{\infty}$
is a frame of subspaces for $\H$.  Also observe that $\span_{i=1,\ldots,\infty}\{W_i^1 , W_i^2 \}
= \H$.  Since for all positive integers $i$ we have
\[{\pi}_{W_i^1}(e_{2i+1}) = \frac{1}{i\sqrt{1+\tfrac{1}{i^2}}}(e_{2i}+
\tfrac{1}{i}e_{2i+1})\quad \mbox{and} \quad
{\pi}_{W_i^2}(e_{2i+1}) = 0,\]
it follows that $\{W_i^1 ,W_i^2 \}_{i=1}^{\infty}$
is not a frame of subspaces for its closed linear span.
\end{example}

Using a Riesz frame of subspaces and Riesz frames for the single
subspaces, we can construct a Riesz frame for $\H$ by just taking
all elements of the Riesz frames.

\begin{prop}
Let $\{W_{i}\}_{i\in I}$ be a Riesz frame of subspaces with respect to
$\{v_i\}_{i \in I}$ for $\H$, and let $\{f_{ij}\}_{j\in I_{i}}$ be a Riesz frame for
$W_i$ with Riesz frame bounds $A$ and $B$ for all $i\in I$. Then
$\{v_i f_{ij}\}_{i\in I,j\in I_i }$ is a Riesz frame for $\H$.
Also, for any $J\subset I$, $\{W_{j}\}_{j\in J}$ is a Riesz frame
of subspaces with respect to $\{v_i\}_{i \in I}$ for its closed linear
span.
\end{prop}

\begin{proof}
Let $C$ and $D$ be the Riesz frame of subspaces bounds for $\{W_{i}\}_{i\in I}$.
For every $i\in I$ choose $J_{i}\subset I_{i}$ and define $\widetilde{W}_i$ by
\[\widetilde{W}_{i} = \span_{j\in J_i}\{f_{ij}\} \subset W_i.\]
Let $f\in \span_{i \in I, j\in J_i}\{f_{ij}\}$.  Then we have
\[\sum_{i\in I, j\in J_i}|\langle f,v_i f_{ij}\rangle |^2 =
\sum_{i\in I}v_i^2 \sum_{j\in J_i}|\langle {\pi}_{\widetilde{W}_i}(f),f_{ij}\rangle |^2
\le \sum_{i\in I}B v_i^2 \|{\pi}_{\widetilde{W}_i}(f)\|^2 \le BD\|f\|^2.\]
Concerning the lower bounds, we compute
\[\sum_{i\in I, j\in J_i}|\langle f,v_i f_{ij}\rangle |^2 =
\sum_{i\in I}v_i^2\sum_{j\in J_i}|\langle {\pi}_{\widetilde{W}_{i}}(f),f_{ij}\rangle |^2
\ge \sum_{i\in I}Av_i^2\|{\pi}_{\widetilde{W}_{i}}(f)\|^2 \ge CA\|f\|^2.\]
\end{proof}


\section{Harmonic frames of subspaces}
\label{harmonic_fos}

Harmonic frames of subspaces are a special case of frames of subspaces, which are equipped
with a natural structure and which occur in several situations, e.g., in Gabor analysis
or in wavelet analysis in the form of multiresolution analysis.

\subsection{The finite case}

We start by giving the definition of a harmonic frame of subspaces for a finite
family of subspaces.

\begin{definition}
A frame of subspaces $\{W_i\}_{i \in I}$ with respect to $\{v_i\}_{i \in I}$ for $\calH$ is a
{\em finite harmonic frame of subspaces with respect to $\{v_i\}_{i \in I}$}, if $I = \{0,\ldots,N-1\}$,
$N \in \NN$ and there exists a unitary operator $U$ on $\calH$ so that
\[U W_{N-1} = W_0 \quad \mbox{and} \quad U W_i = W_{i+1}  \quad \mbox{for all } 0 \le i \le N-2.\]
\end{definition}

If $\{W_i\}_{i \in I}$ is a uniform Parseval frame of subspaces, $U W_{N-1} = W_0$
follows automatically as the following proposition shows. This result equals
the corresponding result in the frame situation \cite[Theorem 4.1]{CK01}.

\begin{theorem}
Let $\{W_i\}_{i \in I}$, $|I| = N$, be a uniform Parseval frame of subspaces
such that there exists a unitary operator $U$ on $\calH$ so that
\[ U W_i = W_{i+1}  \quad \mbox{for all } 0 \le i \le N-2.\]
Then
\[ U W_{N-1} = W_0.\]
\end{theorem}

\begin{proof}
Without loss of generality we can assume that $v_i = 1$ for all $i \in I$.
Let $\{g_j\}_{j \in J}$ be an orthonormal basis for $W_0$. By hypothesis,
$\{U^i g_j\}_{j \in J}$ is an orthonormal basis for $W_{i}$ for all $0 \le i \le N-1$.
Now let $f \in \calH$.
By Proposition \ref{prop_orthonormalbasis} and since $\{W_i\}_{i \in I}$ was assumed to be uniform, we have
\begin{equation} \label{harmonic_eq}
f = \sum_{i =0}^{N-1} \sum_{j \in J} \ip{f}{U^i g_j}U^i g_j.
\end{equation}
Applying $U$, this leads to
\begin{eqnarray*}
Uf & = & \sum_{i =0}^{N-1} \sum_{j \in J} \ip{Uf}{U^i g_j}U^i g_j\\
& = & \sum_{i =0}^{N-1} \sum_{j \in J} \ip{f}{U^{i-1} g_j}U^i g_j\\
& = & U\left[\sum_{i =0}^{N-1} \sum_{j \in J} \ip{f}{U^{i-1} g_j}U^{i-1} g_j\right]\\
& = & U\left[\sum_{i =-1}^{N-2} \sum_{j \in J} \ip{f}{U^i g_j}U^i g_j\right].
\end{eqnarray*}
Since $U$ is unitary and by (\ref{harmonic_eq}), we obtain
\[ f = \sum_{i =-1}^{N-2} \sum_{j \in J} \ip{f}{U^i g_j}U^i g_j
= \sum_{i =0}^{N-1} \sum_{j \in J} \ip{f}{U^i g_j}U^i g_j.\]
This implies
\[ \sum_{j \in J} \ip{f}{U^{-1} g_j}U^{-1} g_j = \sum_{j \in J} \ip{f}{U^{N-1} g_j}U^{N-1} g_j.\]
Now we apply $U$, which yields
\[ \sum_{j \in J} \ip{f}{U^{-1} g_j}g_j = \sum_{j \in J} \ip{f}{U^{N-1} g_j}U^N g_j.\]
Using $U^{-1} f$ instead of $f$ gives
\[ \sum_{j \in J} \ip{f}{g_j}g_j = \sum_{j \in J} \ip{f}{U^N g_j}U^N g_j,\]
which shows that $\pi_{W_0} = \pi_{\span_{j \in J}  \{U^N g_j\}} = \pi_{U W_{N-1}}$. This completes
the proof.
\end{proof}

The following result is \cite[Theorem 4.2 and Theorem 4.3]{CK01}, which we add for completeness in
a reformulated version and with a proof using our results.

\begin{prop}
Let $\varphi \in \H$ and let $V$ be a unitary operator on $\H$ such that
$\{V^j \varphi\}_{j=0}^{M-1}$ is a uniform Parseval frame sequence. Define
the subspace $W_0$ by $W_0 = \ospan_{j=0,\ldots,M-1}\{V^j \varphi\}$. Further let $U$ be
 a unitary operator on $\H$. Then the following conditions are equivalent.
\begin{enumerate}
\item $\{U^i V^j \varphi\}_{i=0,\;\; j=0}^{L-1,\, M-1}$ is a uniform Parseval frame
for $\H$.
\item $\{U^i W_0\}_{i=0}^{L-1}$ is a uniform Parseval
frame of subspaces for $\H$.
\end{enumerate}
\end{prop}

\begin{proof}
This follows immediately from Theorem \ref{transfer_frame} with setting $W_i = U^i W_0$ for all
$i = 1, \ldots, L-1$.
\end{proof}

We conclude this subsection by giving an example for this type of frames of
subspaces.

\begin{example}
Let $a \in \RR$ and $g \in L^2(\RR)$ be such that
$\{E_{am}T_ng\}_{m,n \in \mathbb{Z}}$ is a Parseval Gabor frame. Fix some $N \in \NN$ and
define $W_i$, $0 \le i \le N-1$, by
\[ W_i = \span_{m,n \in \mathbb{Z}} \{E_{a(mN+i)} T_{n} g\}.\]
Then $\{W_i\}_{i=0}^{N-1}$ is a finite harmonic frame of subspaces. To see this let $U := E_a$. Then
\[U W_{N-1} = W_0 \quad \mbox{and} \quad U W_i = W_{i+1}  \quad \mbox{for all } 0 \le i \le N-2.\]
Also,
\[ U^i(E_{amN} T_{n} g) = E_{a(mN+i)} T_{n} g.\]
Hence the sequences $\{E_{a(mN+i)} T_{n} g\}_{m,n \in \mathbb{Z}}$ are even unitarily equivalent
to each other.
\end{example}

Notice that this construction generalizes the one in Example \ref{exa_Gabor}.


\subsection{The infinite case}

We can also define a harmonic frame of subspaces for an infinite family of subspaces.

\begin{definition}
A frame of subspaces $\{W_i\}_{i \in I}$ with respect to $\{v_i\}_{i \in I}$ for $\calH$ is an
{\em infinite harmonic frame of subspaces with respect to $\{v_i\}_{i \in I}$}, if $I = \ZZ$ and
there exists a unitary operator $U$ on $\calH$ so that
\[ U W_i = W_{i+1} \quad \mbox{for all } i \in I.\]
\end{definition}

An interesting example for this type of frame of subspaces is the generalized
frame multiresolution analysis in the sense of Papadakis \cite{Pap03}, which
generalizes the classical multiresolution analysis.

\begin{example}
Here we consider the generalized frame multiresolution analysis in the sense
of Papadakis \cite{Pap03}, whose approach includes all classical multiresolution analysis
(MRAs) in one and higher dimensions as well as the FMRAs of Benedetto and Li \cite{BL98}.

Let $\calH$ be a Hilbert space, $U : \calH \to \calH$ be a unitary operator, and let $G$ be
a unitary abelian group acting on $\calH$. A sequence $\{V_i\}_{i \in \ZZ}$ of closed subspaces
of $\calH$ is a {\em generalized frame multiresolution analysis} of $\calH$, if it satisfies
the following properties.
\renewcommand{\labelenumi}{{\rm (\arabic{enumi})}}
\begin{enumerate}
\item $V_i \subseteq V_{i + 1}$ for all $i \in \ZZ$.
\item $V_i = U^i(V_0)$ for all $i \in \ZZ$.
\item $\bigcap_{i} V_i = \{0\}$ and $\overline{\bigcup_{i} V_i} = \calH$.
\item There exists a countable subset $B$ of $V_0$ such that the set $G(B)=\{g\phi :
g \in G, \phi \in B\}$ is a frame for $V_0$.
\end{enumerate}
\renewcommand{\labelenumi}{{\rm (\roman{enumi})}}
Now $\{W_i\}_{i \in \ZZ}$ is defined by $W_0 := V_1 \cap V_0^\perp$ and
$W_i := U^i(W_0)$ for every $i \in \ZZ$. Then we have
\[ W_i = U^i(V_1 \cap V_0^\perp) = U^i(V_1) \cap U^i(V_0)^\perp = V_{i+1} \cap V_i^\perp.\]

In the case $\calH = L^2(\RR)$, $G = \{T_k : k \in \ZZ\}$, $B$ containing only one element, and
$U = D_2$ being the dilation operator $D_2f(t) = \sqrt{2}f(2t)$ this definition reduces to
the well-known MRA.

In the general situation we have the following relations:
\begin{enumerate}
\item If $\{V_i\}_{i \in \ZZ}$ is a GFMRA of some Hilbert space $\calH$, then $\{W_i\}_{i \in \ZZ}$
is an infinite $1$-uniform harmonic Parseval frame of subspaces, since
$\calH =  \bigoplus_{i\in \ZZ} W_i $, with
unitary operator $U$.
\item Let $\{W_i\}_{i \in \ZZ}$ be an infinite $1$-uniform harmonic Parseval frame of subspaces with unitary
operator denoted by $U$. By Proposition \ref{prop_orthonormalbasis}, we have
$\calH = \bigoplus_{i\in \ZZ} W_i $.
Then we can define $\{V_i\}_{i \in \ZZ}$ by
\[ V_i = \bigoplus_{m \le  i-1} W_m.\]
Then (1)-(3) are obviously satisfied. Therefore $\{W_i\}_{i \in \ZZ}$ is  a GFMRA if and only
if (4) is satisfied.
\end{enumerate}
\end{example}


\section*{Acknowledgments}
The majority of the research for this paper was performed
while the second author was visiting the Departments of Mathematics at
the University of Missouri and Washington University in St.~Louis.
This author thanks these departments for their hospitality and support
during these visits.

The second author is indebted to Damir Baki\'{c}, Ilya Krishtal, Demetrio Labate, Guido Weiss,
and Ed Wilson for helpful discussions.

\bibliographystyle{amsalpha}

\end{document}